\documentclass[12pt]{article}

\usepackage{amssymb,amsmath,amsthm, amsfonts}

=0pt \headheight     = 0pt \headsep        = 0pt \topskip
=0pt
\def\sqr#1#2{{\vcenter{\vbox{\hrule height.#2pt
              \hbox{\vrule width.#2pt height#1pt \kern#1pt \vrule width.#2pt}
              \hrule height.#2pt}}}}
\def\signed #1{{\unskip\nobreak\hfil\penalty50
              \hskip2em\hbox{}\nobreak\hfil#1
              \parfillskip=0pt \finalhyphendemerits=0 \par}}
\def\endpf{\signed {$\sqr69$}}
\def\3n{\negthinspace \negthinspace \negthinspace }
\def\2n{\negthinspace \negthinspace }
\def\1n{\negthinspace }

\def\ds{\displaystyle}

\def\={\buildrel \triangle \over =}

\def\b{\beta}

\def\d{\delta}

\def\l{\lambda}

 \def\n{\nabla}

\def\t{\times}

\def\th{\theta}
\def\o{\omega}

\def\i{\infty}
\def\ns{\noalign{\ss} }
\def\nm{\noalign{\ms} }
\def\G{\Gamma}
\def\D{\Delta}
\def\Th{\Theta}

\def\Si{\Sigma}

\def\O{\Omega}

\def\cF{{\cal F}}

\def\cl{{\cal l}}
\def\no{\noindent}

\def\ss{\smallskip}
\def\ms{\medskip}
\def\bs{\bigskip}
\def\q{\quad}
\def\qq{\qquad}
\def\hb{\hbox}

\def\lan{\mathop{\langle}}
\def\ran{\mathop{\rangle}}

\def\max{\mathop{\rm max}}
\def\min{\mathop{\rm min}}
\def\exp{\mathop{\rm exp}}

\def\pa{\partial}

\def\cd{\cdot}

\def\ae{\hbox{\rm a.e.{ }}}
\def\as{\hbox{\rm a.s.{ }}}

\def\cl{\overline}

\def\|{\Big |}
\def\({\Big (}
\def\){\Big )}
\def\[{\Big[}
\def\]{\Big]}
\def\be{\begin{equation}}
\def\bel{\begin{equation}\label}
\def\ee{\end{equation}}
\def\bt{\begin{theorem}}
\def\bcd{\begin{condition}}
\def\ecd{\end{condition}}
\def\et{\end{theorem}}
\def\bc{\begin{corollary}}
\def\ec{\end{corollary}}
\def\bde{\begin{definition}}
\def\ede{\end{definition}}
\def\bl{\begin{lemma}}
\def\el{\end{lemma}}
\def\bp{\begin{proposition}}
\def\ep{\end{proposition}}
\def\br{\begin{remark}}
\def\er{\end{remark}}
\def\ba{\begin{array}}
\def\ea{\end{array}}
\def\ed{\end{document}}
\def\ns{\noalign{\ms}}
\def\ds{\displaystyle}

\def\square#1{\vbox{\hrule\hbox{\vrule height#1%
     \kern#1\vrule}\hrule}}
\def\rectangle#1#2{\vbox{\hrule\hbox{\vrule height#1%
     \kern#2\vrule}\hrule}}

\font\tenbb=msbm10 \font\sevenbb=msbm7 \font\fivebb=msbm5

\newfam\bbfam
\scriptscriptfont\bbfam=\fivebb \textfont\bbfam=\tenbb
\scriptfont\bbfam=\sevenbb

\newtheorem{lemma}{Lemma}[section]
\newtheorem{remark}{Remark}[section]

\newtheorem{theorem}{Theorem}[section]
\newtheorem{corollary}{Corollary}[section]

\newtheorem{definition}{Definition}[section]
\newtheorem{proposition}{Proposition}[section]
\newtheorem{condition}{Condition}[section]

\makeatletter
   
   \@addtoreset{equation}{section}
\makeatother

\begin{document}
\title{\bf Carleman and Observability Estimates for Stochastic Wave Equations\thanks{This work was
supported by the NSF of China under grant 10525105 and the NCET of
China under grant NCET-04-0882. Part of this work was done when the
author visited the Shanghai Key Laboratory for Contemporary Applied
Mathematics at Fudan University. The author acknowledges Professor
Shanjian~Tang for stimulating discussion. \ms}}

\author{Xu Zhang\thanks{Key Laboratory of Systems and Control, Academy
of Mathematics and Systems Sciences, Chinese Academy of Sciences,
Beijing 100080, China; and Yangtze Center of Mathematics, Sichuan
University, Chengdu 610064, China. {\small\it E-mail:} {\small\tt
xuzhang@amss.ac.cn}.} }

\date{}

\maketitle

\begin{abstract}
\no Based on a fundamental identity for stochastic hyperbolic-like
operators, we derive in this paper a global Carleman estimate (with
singular weight function) for stochastic wave equations. This leads
to an observability estimate for stochastic wave equations with
non-smooth lower order terms. Moreover, the observability constant
is estimated by an explicit function of the norm of the involved
coefficients in the equation.
\end{abstract}

\bs

\no{\bf 2000 Mathematics Subject Classification}.  Primary 35B45;
Secondary 93B07, 35R30.

\bs

\no{\bf Key Words}.  Carleman estimate, singular weight function,
observability estimate, stochastic wave equation.


\section{Introduction and main results}

Let $T>0$, $G \subset \mathbb{R}^{n}$ ($n\in\mathbb{N}$) be a given
bounded domain with a $C^2$ boundary $\G$. Fix any
$x_0\in\mathbb{R}^d\setminus \cl G$.  It is clear that
 \bel{4.1} 0<R_0\=\min_{x\in G}|x-x_0|<
R_1\=\max_{x\in G}|x-x_0|.
 \ee
Put
 \be\label{boundary}
 \G_0\=\big\{x\in \G\;\big|\; (x-x_0)\cdot\nu(x)>0\big\},
 \ee
where $\nu(x)$ is the unit outward normal vector of $G$ at
 $x\in \G$. Also, put $ Q\=(0,T)\times G$, $ \Si\=(0,T)\times \G$ and $ \Si_0\=(0,T)\times \G_0$ . Throughout
this paper, we will use $C$ to denote a generic positive constant
depending only on $T$, $G$ and $G_0$, which may change from line to
line.

Let $(\O,\cF,\{\cF_t\}_{t\ge 0},P)$  be a complete filtered
probability space on which a one dimensional standard Brownian
motion $\{w(t)\}_{t\ge0}$ is defined. Let $H$ be a Banach space. We
denote by $L_{\cF}^2(0,T;H)$ the Banach space consisting of all
$H$-valued $\{\cF_t\}_{t\ge 0}$-adopted processes $X(\cd)$ such that
$\mathbb{E}(|X(\cd)|_{L^2(0,T;H)}^2)<\i$, with the canonical norm;
by$L_{\cF}^\i(0,T;H)$ the Banach space consisting of all $H$-valued
$\{\cF_t\}_{t\ge 0}$-adopted bounded processes; and by
$L_{\cF}^2(\O;C([0,T];H))$ the Banach space consisting of all
$H$-valued $\{\cF_t\}_{t\ge 0}$-adopted continuous processes
$X(\cd)$ such that $\mathbb{E}(|X(\cd)|_{C([0,T];H)}^2)<\i$, with
the canonical norm.

Assume
 \bel{wj1}
 \ba{ll}
 a_1\in L^{\i}_{\cF}(0,T;L^{\i}(G)),\qq& a_2\in
 L^{\i}_{\cF}(0,T;L^{\i}(G;\mathbb{R}^n)), \\
 \ns
 a_3\in L^{\i}_{\cF}(0,T;L^n(G)),\qq &a_4\in
 L^{\i}_{\cF}(0,T;L^\i(G)),
 \ea
 \ee
and
 \bel{wj2}
 f\in L^2_{\cF}(0,T;L^2(G)),\qq g\in L^2_{\cF}(0,T;L^2(G)).
 \ee
Let us consider the following stochastic wave equation:
 \bel{hh6.1}
 \left\{
 \ba{ll}
 \ds dy_t-\D ydt=(a_1y_t+\lan a_2,\n y\ran+a_3y+f)dt+(a_4y+g)
 dw(t)&\hb{ in }Q,\\
 \ns
 y=0&\hb{ on }\Si,\\
 \ns
 y(0)=y_0,\q y_t(0)=y_1&\hb{ in }G.
 \ea
 \right.
 \ee
Here, we denote the scalar product in $\mathbb{R}^n$ by
$\lan\cd,\cd\ran$. For any initial data
 \bel{w2}
 (y_0,y_1)\in L^2(\O,\cF_0,P;H_0^1(G)\times L^2(G)),
 \ee
it is easy to show that system (\ref{hh6.1}) admits one and only one
weak solution
 $$y\in
L_{\cF}^2(\O;C([0,T];H_0^1(G))\bigcap C^1([0,T];L^2(G))).
 $$
By means of the classical multiplier approach and energy estimate,
following \cite{Lions, Zh4}, it is not difficult to show the
following hidden regularity for the solution of system (\ref{hh6.1})
(Hence we omit the details):

\bp\label{p1}
Under assumptions (\ref{wj1}), (\ref{wj2}) and (\ref{w2}), the
solution of system (\ref{hh6.1}) satisfies $\frac{\pa y}{\pa\nu}\in
L^2_{\cF}(0,T;L^2(\G))$. Moreover
 \bel{w4}
 \ba{ll}\ds
 \left|\frac{\pa y}{\pa\nu}\right|_{ L^2_{\cF}(0,T;L^2(\G))}\\
 \ns
 \ds\le C\[|(y_0,y_1)|_{L^2(\O,\cF_0,P;H_0^1(G)\times L^2(G))}+|f|_{L^2_{\cF}(0,T;L^2(G)}+|g|_{
L^2_{\cF}(0,T;L^2(G))}\]\\
\ns
\ds\q\times\exp\Big\{C\[|(a_1,a_4)|_{L^{\i}_{\cF}(0,T;(L^\i(G))^2)}^2+|
a_2|_{
 L^{\i}_{\cF}(0,T;L^{\i}(G;\mathbb{R}^n))}^2+|a_3|_{L^{\i}_{\cF}(0,T;L^n(G)))}^2\]\Big\}. \ea
 \ee
\ep

The main purpose of this paper is to derive a boundary observability
estimate for system (\ref{hh6.1}). For this, we choose a
sufficiently small constant $c\in (0,1)$ so that (Recall (\ref{4.1})
for $R_0$ and $R_1$)
 $$
 \frac{(4+5c)R_0^2}{9c}>R_1^2.
 $$
Then, in the sequel, we take $T(>2R_1)$ sufficiently large such that
 \bel{2z}
\frac{4(4+5c)R_0^2}{9c}>c^2T^2>4R_1^2.
 \ee
Our observability estimate for system (\ref{hh6.1}) is stated as
follows:

\bt\label{t1}
Let (\ref{wj1})--(\ref{wj2}) hold, $R_1$ and $\G_0$ be given
respectively by (\ref{4.1}) and (\ref{boundary}), and $T$ satisfy
(\ref{2z}). Then solutions of system (\ref{hh6.1}) satisfy
 \bel{12e2}
 \ba{ll}\ds
 |(y(T),y_t(T))|_{L^2(\O,\cF_T,P;H_0^1(G)\times L^2(G))}\\
 \ns
 \ds\le C\left[\left|\frac{\pa y}{\pa\nu}\right|_{ L^2_{\cF}(0,T;L^2(\G_0))}+|f|_{L^2_{\cF}(0,T;L^2(G)}+|g|_{
L^2_{\cF}(0,T;L^2(G))}\right]\\
\ns
\ds\q\times\exp\Big\{C\[|(a_1,a_4)|_{L^{\i}_{\cF}(0,T;(L^\i(G))^2)}^2+|
a_2|_{
 L^{\i}_{\cF}(0,T;L^{\i}(G;\mathbb{R}^n))}^2+|a_3|_{L^{\i}_{\cF}(0,T;L^n(G))}^2\]\Big\},\\
 \ns
 \ds\qq\qq\qq\qq\forall\;(y_0,y_1)\in L^2(\O,\cF_0,P;H_0^1(G)\times L^2(G)).
 \ea
 \ee
\et

It is well-known that observability estimate is an important tool
for the study of stabilization and controllability problems for
deterministic PDEs. We refer to \cite{Zu} for a recent survey in
this respect. Although there are numerous references addressed to
the observability problems for deterministic PDEs, very little is
known for the stochastic counterpart and it remains to be further
understood. Indeed, to the best of our knowledge, \cite{BRT} is the
only one publication in this field, which is devoted to the
controllability/observability for the stochastic heat equation. As
far as we know, nothing is known for the observability estimate on
the stochastic wave equation.

Similar to the deterministic setting, we shall use a stochastic
version of the global Carleman estimate to establish inequality
(\ref{12e2}). The difficulty to do this is the very fact that,
unlike the deterministic situation, system (\ref{hh6.1}), a
stochastic wave equation, is {\it time-irreversible}. Therefore, one
can not simply mimic the usual Carleman inequality for the
deterministic wave equations (See \cite{FYZ, Zh4} and the references
cited therein). Rather, instead of the usual smooth weight function,
one has to introduce another singular weight function to derive the
desired Carleman estimate for system (\ref{hh6.1}).

More precisely, for any (large) $\l>0$ and any (small) $c>0$, set
 \bel{a3'}
\ell=\ell(t,x)\=\l\[|x-x_0|^2-c\(t-{T\over2}\)^2\], \qq \th\=e^\ell.
 \ee
Also, for any $\b>0$, we set
 \bel{4.10}
 \Th=\Th(t)\=\exp\left\{-\frac{\b }{t(T-t)}\right\},\qq 0<t<T.
 \ee
It is easy to see that $\Th(t)$ decays rapidly to $0$ as $t\to0$ or
$t\to T$. Our Carleman estimate for system (\ref{hh6.1}) is stated
as follows:

\bt\label{t2}
Let (\ref{wj1})--(\ref{wj2}) hold, $R_1$ and $\G_0$ be given
respectively by (\ref{4.1}) and (\ref{boundary}), and $T$ satisfy
(\ref{2z}).  Then there exist a constant $\b>0$ (which is very
small), and a constant
$$\l^*=C\[1+|(a_1,a_4)|_{L^{\i}_{\cF}(0,T;(L^\i(G))^2)}^2+|
a_2|_{
 L^{\i}_{\cF}(0,T;L^{\i}(G;\mathbb{R}^n))}^2+|a_3|_{L^{\i}_{\cF}(0,T;L^n(G))}^2\],
 $$
such that solutions of system (\ref{hh6.1}) satisfy
 \bel{12e22}
 \ba{ll}
\ds \l\mathbb{E}\int_Q\Th \th^2(y_t^2+|\n y|^2+\l^2
y^2)dxdt\\
 \ns
 \le \ds C\mathbb{E}\left\{\l\int_{\Si_0}\Th\th^2\left|\frac{\pa y}{\pa\nu}\right|^2d\Si_0+\int_Q\Th\th^2(f^2+\l g^2)dxdt\right\},\\
 \ns
 \ds\qq\qq\qq\qq\forall\;(y_0,y_1)\in L^2(\O,\cF_0,P;H_0^1(G)\times L^2(G)),\q\forall\;\l\ge \l^*.
 \ea
 \ee
\et

Carleman estimate is a fundamental tool for the study of control and
inverse problems for deterministic PDEs (\cite{I, Zu}). Similar to
the situation for observability estimate, although there are
numerous references addressed to Carleman estimate for deterministic
PDEs, to the best of our knowledge, \cite{BRT, TZ} are the only two
references for the stochastic counterpart, which are devoted to the
stochastic heat equation. It would be quite interesting to extend
the deterministic Carleman estimate for other PDEs to the stochastic
ones, but there are many things to be done, and some of which seem
to be challenging. In this paper, in order to present the key idea
in the simplest way, we do not pursue the full technical generality.

The rest of this paper is organized as follows. In Section \ref{s2},
as a key preliminary, we present an identity for a stochastic
hyperbolic-like operator. Then, in Section \ref{s3}, we derive
pointwise Carleman-type estimates for the stochastic wave operator.
Finally, Section \ref{s4} is devoted to the proof of Theorems
\ref{t1}-\ref{t2}.

\section{Identity for a stochastic hyperbolic-like operator}\label{s2}

For simplicity, we denote $\ds\sum^{n}_{i,j=1}$ and
$\ds\sum^{n}_{i=1}$ simply by $\ds\sum_{i,j}$ and $\ds\sum_{i}\,$,
respectively. Also, we will use the notation $ u_i= u_{ x_i}$, where
$x_i$ is the $i$-th coordinate of a generic point
$x=(x_1,\cdots,x_n)$ in $\mathbb{R}^n$. In a similar manner, we use
the notation $\ell_i$, $v_i$, etc. for the partial derivatives of
$\ell$ and $v$ with respect to $x_i$.

We show the following fundamental identity for a stochastic
hyperbolic-like operator:

\bt\label{c1t1}
Let $b^{ij}\in C^1((0,T)\t\mathbb{R}^n)$ satisfying
 \bel{1}
 b^{ij}=b^{ji},\qq i,j=1,2,\cdots,n,
 \ee
$u,\ \ell,\ \Psi\in C^2((0,T)\t\mathbb{R}^n)$. Assume $u$ is a
$H^2_{loc}(\mathbb{R}^n)$-valued $\{\cF_t\}_{t\ge 0}$-adopted
processes such that $u_t$ is a $L^2_{loc}(\mathbb{R}^n)$-valued
semi-martingale.
Set
 $\th=e^{\ell }$ and $v=\th u$. Then for $\ae$ $x\in \mathbb{R}^n$ and $P$-$\as$ $\o\in\O$,
 \bel{c1e2a}
 \3n\3n\ba{ll}
 \displaystyle
 \th\(-2\ell_tv_t+2\sum_{i,j}b^{ij}\ell_iv_j+\Psi v\)\[du_t-\sum_{i,j}(b^{ij}u_i)_jdt\]\\
 \noalign{\ss}
 \displaystyle\q+
  \sum_{i,j}\[\sum_{i',j'}\(2b^{ij} b^{i'j'}\ell_{i'}v_iv_{j'}
 -b^{ij}b^{i'j'}\ell_iv_{i'}v_{j'}\)-2b^{ij}\ell_tv_iv_t
 +b^{ij}\ell_iv_t^2\\
 \noalign{\ss}
 \displaystyle\qq\qq\qq+\Psi b^{ij}v_iv-\(A\ell_i+{\Psi_i\over2}\)b^{ij} v^2\]_jdt\\
 \noalign{\ss}
 \displaystyle\q+d\[\sum_{i,j}
 b^{ij}\ell_tv_iv_j-2\sum_{i,j}b^{ij}\ell_iv_jv_t+\ell_tv_t^2-\Psi
 v_tv+\(A\ell_t+{\Psi_t\over2}\)v^2\]\\
 \noalign{\ss}
 \displaystyle
 =\Big\{\[\ell_{tt}+\sum_{i,j}(b^{ij}
 \ell_i)_j-\Psi\]v_t^2-2\sum_{i,j}\[(b^{ij}\ell_j)_t+b^{ij}\ell_{tj}\]v_iv_t\\
 \noalign{\ss}
 \displaystyle\q+\sum_{i,j} \Big\{(b^{ij}\ell_t)_t+\sum_{i',j'}\[2b^{ij'}(b^{i'j}\ell_{i'})_{j'} -
 (b^{ij}b^{i'j'}\ell_{i'})_{j'}\]+\Psi b^{ij}
 \Big\}v_iv_j\\
 \noalign{\ss}
 \displaystyle\q
 +Bv^2+\(-2\ell_tv_t+2\sum_{i,j}b^{ij}\ell_iv_j+\Psi v\)^2\Big\}dt+\th^2\ell_t(du_t)^2,
 \ea
 \ee
where
\bel{c1e3a}
\left\{
 \ba{ll}
 \ds A\=(\ell_t^2-\ell_{tt})-\sum_{i,j}
 (b^{ij}\ell_i\ell_j-b^{ij}_j\ell_i
 -b^{ij}\ell_{ij})-\Psi,\\
  \ns
 \ds
 B\=A\Psi+(A\ell_t)_t-
 \sum_{i,j}(Ab^{ij}\ell_i)_j +{1\over
 2}\[\Psi_{tt}-\sum_{i,j} (b^{ij}\Psi_i)_j\].
  \ea
\right.
\ee
\et

{\it Proof.}  Recall that
 $$
 v(t,x)=\th(t,x)u(t,x).
 $$
Hence $u_t=\th^{-1}(v_t-\ell _tv)$ and $u_j=\th^{-1}(v_j-\ell
_jv)$ for $j=1,2,\cdots,n$. Hence,
 \bel{2c2t1}
 \ba{ll}
 \displaystyle
 du_{t}=\th^{-1}[d v_{t}-2\ell _tv_tdt+(\ell _t^2-\ell _{tt})vdt].
 \ea
 \ee
Similarly, by symmetry condition (\ref{1}), one may check that
 \bel{w2.1}
 \displaystyle
 \sum_{i,j}(b^{ij}u_i)_j
 =\th^{-1}\sum_{i,j}\[(b^{ij}v_i)_j-2b^{ij}\ell_iv_j
 +(b^{ij}\ell_i\ell_j-b^{ij}_j\ell_i
 -b^{ij}\ell_{ij})v\].
 \ee
Therefore, by (\ref{2c2t1})--(\ref{w2.1}), and recalling the
definition of $A$ in (\ref{c1e3a}), we get
 \bel{c1e5}
 \ba{ll}
 \displaystyle\th\(-2\ell_tv_t+2\sum_{i,j}b^{ij}\ell_iv_j+\Psi v\)\[du_t-\sum_{i,j}(b^{ij}u_i)_jdt\]\\
   \nm
 \ds= \(-2\ell_tv_t+2\sum_{i,j}b^{ij}\ell_iv_j+\Psi v\)\Big\{dv_t-\[\sum_{i,j} (b^{ij}v_i)_j-Av\\
  \nm
 \ds\qq\q+2\ell_tv_t-2\sum_{i,j}b^{ij}\ell_iv_j-\Psi v\]dt\Big\}\\
   \nm
 \ds= \(-2\ell_tv_t+2\sum_{i,j}b^{ij}\ell_iv_j+\Psi v\)dv_t\\
   \nm
 \ds\q
 +\(-2\ell_tv_t+2\sum_{i,j}b^{ij}\ell_iv_j+\Psi v\)\[-\sum_{i,j} (b^{ij}v_i)_j+Av\]dt\\
  \nm
 \ds\q+\(-2\ell_tv_t+2\sum_{i,j}b^{ij}\ell_iv_j+\Psi v\)^2dt.
 \ea
 \ee

We now analyze the first two terms in the right hand side of
(\ref{c1e5}).

First, using It\^o's formula, we have
 \bel{w1}
 \ba{ll}
 \ds\(-2\ell_tv_t+2\sum_{i,j}b^{ij}\ell_iv_j+\Psi v\)dv_t\\
 \ns
 \ds=d\[ \(-\ell_tv_t+2\sum_{i,j}b^{ij}\ell_iv_j+\Psi
 v\)v_t\]\\
 \ns
 \ds\q-\[-\ell_{tt}v_t^2+2\sum_{i,j}(b^{ij}\ell_i)_tv_jv_t+2\sum_{i,j}b^{ij}\ell_iv_{tj}v_t+\Psi
 v_t^2+\Psi_tvv_t\]dt+\ell_t(dv_t)^2\\
 \ns
 \ds= d\(-\ell_tv_t^2+2\sum_{i,j}b^{ij}\ell_iv_jv_t+\Psi
 vv_t-{\Psi_t\over2}v^2\)\\
 \ns
 \ds\q+\Big\{- \sum_{i,j}(b^{ij}\ell_iv_t^2)_j+ \[\ell_{tt}+\sum_{i,j}(b^{ij}
 \ell_i)_j-\Psi\]v_t^2-2\sum_{i,j}(b^{ij}\ell_j)_tv_iv_t+{\Psi_{tt}\over
 2}v^2\Big\}dt\\
 \ns
 \ds\q+ \th^2\ell_t(du_t)^2.
 \ea
 \ee

Next,
 \bel{1c1.6}
 \ba{ll}
 \displaystyle -2\ell _tv_t\[-\sum_{i,j}
(b^{ij}v_i)_j+Av\]\\
 \noalign{\ss}
 \displaystyle=2\[\sum_{i,j}
 (b^{ij}\ell_tv_iv_t)_j-\sum_{i,j}
 b^{ij}\ell_{tj}v_iv_t\]-\sum_{i,j}
 b^{ij}\ell_t(v_iv_j)_t-A\ell_t(v^2)_t\\
 \noalign{\ss}
 \displaystyle=2\[\sum_{i,j}
 (b^{ij}\ell_tv_iv_t)_j-\sum_{i,j}
 b^{ij}\ell_{tj}v_iv_t\]+\sum_{i,j}
 (b^{ij}\ell_t)_tv_iv_j\\
 \noalign{\ss}
 \displaystyle\q-\(\sum_{i,j}
 b^{ij}\ell_tv_iv_j+A\ell_tv^2\)_t+(A\ell_t)_tv^2.
 \ea
 \ee

Further, by means of a direct computation, one may check that
 \bel{1c1.3}
 \3n\1n\ba{ll}
 \displaystyle
 2\sum_{i,j}b^{ij}\ell_iv_j\[-\sum_{i,j} (b^{ij}v_i)_j+Av\]\\
 \noalign{\ss}
 \displaystyle
 = - \sum_{i,j}\[\sum_{i',j'}\(2b^{ij} b^{i'j'}\ell_{i'}v_iv_{j'}
 -b^{ij}b^{i'j'}\ell_iv_{i'}v_{j'}\)-Ab^{ij}\ell_i v^2\]_j\\
 \noalign{\ss}
 \displaystyle
 \q+\sum_{i,j,i',j'} \[2b^{ij'}(b^{i'j}\ell_{i'})_{j'} -
 (b^{ij}b^{i'j'}\ell_{i'})_{j'}\]v_iv_j-
 \sum_{i,j}(Ab^{ij}\ell_i)_j v^2,
 \ea
 \ee
and
 \bel{2c2t11}
 \ba{ll}
 \displaystyle \Psi v\[-\sum_{i,j}
 (b^{ij}v_i)_j+Av\]&\ds= -\sum_{i,j}
 \(\Psi b^{ij}v_iv-{\Psi_i\over2}b^{ij} v^2\)_j+\Psi \sum_{i,j}
 b^{ij}v_iv_j\\
 \noalign{\ss}
 &\displaystyle\q+\[-{1\over2}\sum_{i,j} (b^{ij}\Psi_i)_j+A\Psi\] v^2.
 \ea
 \ee

Finally, combining (\ref{c1e5})--(\ref{2c2t11}), we arrive at the
desired equality (\ref{c1e2a}). \endpf

\section{Pointwise Carleman-type estimates for the stochastic wave
operator}\label{s3}

In this section, we show a pointwise Carleman-type estimate (with
singular weight) for the stochastic wave operator ``$du_t-\D udt$".

To begin with, by taking $(b^{ij})_{n\t n}=I$, the identity matrix,
and $\th=e^\ell$ (with $\ell$ given in (\ref{a3'})) in Theorem
\ref{c1t1}, one has the following pointwise Carleman-type estimate
for the stochastic wave operator.
 \bl\label{l1}
Let $u,\ \ell,\ \Psi\in C^2((0,T)\t\mathbb{R}^n)$ and
$k\in\mathbb{R}$. Assume $u$ is a $H^2_{loc}(\mathbb{R}^n)$-valued
$\{\cF_t\}_{t\ge 0}$-adopted processes such that $u_t$ is a
$L^2_{loc}(\mathbb{R}^n)$-valued semi-martingale.
Set $v=\th u$. Then for $\ae$ $x\in \mathbb{R}^n$ and $P$-$\as$
$\o\in\O$, it holds
 \bel{b3}
 \ba{ll}
\displaystyle
 \th(-2\ell_tv_t+2\n\ell\cd\n v+\psi v)(du_t-\D udt)\\
 \noalign{\ss}
 \displaystyle\q+d\[\ell_t(v_t^2+|\n v|^2)-2(\n \ell)\cd (\n v)v_t
 -\Psi vv_t+A\ell_tv^2\] \\
\noalign{\ss} \q\displaystyle+\sum_{i=1}^n\Big\{2v_i(\n\ell)\cd(\n
v)
 -\ell_i|\n v|^2-2\ell_tv_tv_i+\ell_iv_t^2+\Psi vv_i
 -A\ell_iv^2\Big\}_idt\\
 \noalign{\ss}
 \displaystyle
\ge \[(1-k)\l v_t^2+(k+3-4c)\l|\n
v|^2+Bv^2\\
 \noalign{\ss}
 \displaystyle\q+\(-2\ell_tv_t+2\n\ell\cd\n v+\psi v\)^2\]dt+\th^2\ell_t(du_t)^2,
 \ea
 \ee
where
 \bel{c4}
\left\{
\ba{ll}
\ds\Psi\=(2n-2c-1+k)\l, \\
 \ns
\ds A=4\left[c^2\left(t-\frac{T}{2}\right)^2-|x-x_0|^2\right]\l^2+\l(4c+1-k),\\
 \ns
\ds
 B=4\left[(4c+5-k)|x-x_0|^2-(8c+1-k)c^2\left(t-\frac{T}{2}\right)^2\right]\l^3+O(\l^2).
 \ea
\right.
\ee
\el

The desired pointwise Carleman-type estimate (with singular weight
function $\Th$) for the stochastic wave operator reads as follows:

\bt\label{4t1}
Let $u\in C^2([0,T]\t\cl\O)$, $v=\th u$, and $T$ satisfy (\ref{2z}).
Then there exist three constant $\l_0>0$, $\b_0>0$ and $c_0>0$,
independent of $u$, such that for all $\b\in (0,\b_0)$ and $\l\ge
\l_0$ it
 holds
 \bel{4.11}
 \ba{ll}
 \ds\Th\th(-2\ell_tv_t+2\n\ell\cd\n v+\psi v)(du_t-\D udt)\\
 \noalign{\ss}
 \displaystyle\q+d\left\{\Th\[\ell_t(v_t^2+|\n v|^2)-2(\n \ell)\cd (\n v)v_t
 -\Psi vv_t+A\ell_tv^2\]\right\} \\
\noalign{\ss}
\q\displaystyle+\sum_{i=1}^n\Big\{\Th\[2v_i(\n\ell)\cd(\n v)
 -\ell_i|\n v|^2-2\ell_tv_tv_i+\ell_iv_t^2+\Psi vv_i
 -A\ell_iv^2\]\Big\}_idt\\
 \noalign{\ss}
 \displaystyle
\ge \left[c_0\l\Th \th^2(u_t^2+|\n u|^2+\l^2
u^2)+\Th\(-2\ell_tv_t+2\n\ell\cd\n v+\psi v\)^2\right]dt+\Th
\th^2\ell_t(du_t)^2,
 \ea\ee
with $A$ and $\Psi$ given by (\ref{c4}).
\et

\br
The main difference between the pointwise estimates (\ref{b3}) and
(\ref{4.11}) is that we introduce a singular ``pointwise" weight in
(\ref{4.11}).  Another difference between (\ref{b3}) and
(\ref{4.11}) is that $T$ is arbitrary in the former estimate; while
for the later one needs to take $T$ to be large enough.
\er

{\it Proof of Theorem \ref{4t1}.} We use some idea in the proof of
\cite[Theorem 1]{ZZ}. The proof is divided it into several steps.

\ms

{\it Step 1.} We multiply both sides of inequality (\ref{b3}) by
$\Th$. Obviously, we have (recall (\ref{c4}) for $A$ and $\Psi$)
 \bel{4.12}
 \ba{ll}
 \displaystyle\Th d\[\ell_t(v_t^2+|\n v|^2)-2(\n \ell)\cd (\n v)v_t
 -\Psi vv_t+A\ell_tv^2\]\\
 \noalign{\ms}\displaystyle= d\left\{\Th\[\ell_t(v_t^2+|\n v|^2)-2(\n \ell)\cd (\n v)v_t
 -\Psi vv_t+A\ell_tv^2\]\right\}\\
 \noalign{\ms}\displaystyle\q-\frac{\b  (T-2t)}{t^2(T-t)^2}\Th\[\ell_t(v_t^2+|\n v|^2)-2(\n \ell)\cd (\n v)v_t
 -\Psi vv_t+A\ell_tv^2\]dt.
 \ea
 \ee
Note that
 \bel{4.13}
 \ba{ll}
 \displaystyle\|-\frac{\b  (T-2t)}{t^2(T-t)^2}\Th\[-2(\n \ell)\cd (\n v)v_t
 -\Psi vv_t\]\|\\
 \noalign{\ms}\displaystyle\le\frac{\b  |T-2t|}{t^2(T-t)^2}\Th\[2 |(\n \ell)\cd (\n v)v_t|
 +|\Psi v_tv|\]\\
 \noalign{\ms}\displaystyle\le \frac{\b  |T-2t|}{t^2(T-t)^2}
 \Th \Big[(|\n \ell|+1)v_t^2
 +|\n \ell||\n v|^2+{1\over{4}}\Psi^2v^2\Big].
 \ea\ee
Thus by (\ref{b3}), and using (\ref{4.12})--(\ref{4.13}), we get
 \bel{4.14}
 \ba{ll}
 \Th\th(-2\ell_tv_t+2\n\ell\cd\n v+\psi v)(du_t-\D udt)\\
 \noalign{\ss}
 \displaystyle\q+d\left\{\Th\[\ell_t(v_t^2+|\n v|^2)-2(\n \ell)\cd (\n v)v_t
 -\Psi vv_t+A\ell_tv^2\]\right\} \\
\noalign{\ss}
\q\displaystyle+\sum_{i=1}^n\Big\{\Th\[2v_i(\n\ell)\cd(\n v)
 -\ell_i|\n v|^2-2\ell_tv_tv_i+\ell_iv_t^2+\Psi vv_i
 -A\ell_iv^2\]\Big\}_idt\\
 \noalign{\ss}
 \displaystyle\ge \left\{\Th (1-k)\l v_t^2+\Th(k+3-4c)\l|\n v|^2+\frac{\b  (T-2t)}{t^2(T-t)^2}\ell _t\Th(v_t^2+|\n v|^2)\right.\\
 \noalign{\ss}
 \displaystyle\q-\frac{\b  |T-2t|}{t^2(T-t)^2}\Th\[(|\n \ell|+1)v_t^2
 +|\n \ell||\n v|^2\]\\
 \noalign{\ss}
 \displaystyle\left.\q+\[B+\frac{\b  (T-2t)}{t^2(T-t)^2}\ell _tA
 -\frac{\b  |T-2t|}{4t^2(T-t)^2}\Psi^2\]\Th v^2+\Th\(-2\ell_tv_t+2\n\ell\cd\n v+\psi v\)^2\right\}dt\\
 \ns
 \q\ds+\Th \th^2\ell_t(du_t)^2,
 \ea
 \ee
where $B$ is given by (\ref{c4}).

\ms

{\it Step 2.} Recalling that $\ell$ and $\Psi$ are given
respectively by (\ref{a3'}) and (\ref{c4}), we get
 \bel{4.15}
 \ba{ll}\ds
 \hb{\rm RHS of }(\ref{4.14})=&\ds\left[\l\Th (F_1v_t^2+F_2|\n v|^2)+\l^3\Th Gv^2+\Th\(-2\ell_tv_t+2\n\ell\cd\n v+\psi v\)^2\right]dt\\
 \ns&\ds+\Th \th^2\ell_t(du_t)^2,
 \ea
 \ee
where
 \bel{4.16}
 F_1\= 1-k+\frac{c\b  (T-2t)^2}{t^2(T-t)^2}
 -\frac{\b  |T-2t|}{t^2(T-t)^2}(2|x-x_0|+\l^{-1}),
 \ee
 \bel{xx4.16}
 F_2\=k+3-4c+\frac{c\b  (T-2t)^2}{t^2(T-t)^2}
 -\frac{2\b
 |T-2t||x-x_0|}{t^2(T-t)^2},\q\;\,
 \ee
and
 \bel{4.17}
 \ba{ll}
 G\=&\ds 4\left[(4c+5-k)|x-x_0|^2-(8c+1-k)c^2\left(t-\frac{T}{2}\right)^2\right]+O(\l^{-1})\\
 \noalign{\ss}
 &\ds\q+\frac{\b
 |T-2t|}{t^2(T-t)^2}
 \Big\{4c|T-2t|\[c^2(t-T/2)^2-|x-x_0|^2\]+O(\l^{-1})\Big\}.
\ea \ee

{\it Step 3.} Let us show that $F_1$, $F_2$ and $G$ are positive
when $\l$ is large enough and $\b$ is sufficiently small. For this,
put
 $$
 \ba{ll}\ds
 F_1^0\=1-k, \qq F_2^0\=k+3-4c, \\
 \ns
 \ds
 G^0\=4\left[(4c+5-k)|x-x_0|^2-(8c+1-k)c^2\left(t-\frac{T}{2}\right)^2\right]+O(\l^{-1}),
 \ea
 $$
which are respectively the nonsingular part of $F_1$, $F_2$ and $G$.
Similarly, put
 $$
 \ba{ll}
 \ds F_1^1\=\ds\frac{c\b  (T-2t)^2}{t^2(T-t)^2}
 -\frac{\b  |T-2t|}{t^2(T-t)^2}(2|x-x_0|+\l^{-1}),\qq F_2^1\=\frac{c\b  (T-2t)^2}{t^2(T-t)^2}
 -\frac{2\b
 |T-2t||x-x_0|}{t^2(T-t)^2},\\
  \ns
  \ds G^1\ds\=\frac{\b
 |T-2t|}{t^2(T-t)^2}
\Big\{4c|T-2t|\[c^2(t-T/2)^2-|x-x_0|^2\]+O(\l^{-1})\Big\},
 \ea
 $$
which are respectively the singular part of $F_1$, $F_2$ and $G$.

Further, we choose $k=1-c$. It is easy to see that both $F_1^0$ and
$F_1^0$ are positive, and
 $$
 G^0\ge 4(4+5c)R_0^2-9c^3T^2+O(\l^{-1}),
 $$
which, via the first inequality in (\ref{2z}), is positive provided
that $\l$ is sufficiently large.

When $t$ is close to $0$ or $T$, i.e., $t\in I_0\=(0,\d_0)\cup
(T-\d_0,T)$ for some sufficiently small $\d_0\in (0,T/2)$,  the
dominant terms in $F_i$ ($i=1,2$) and $G$ are the singular ones. For
$t\in I_0$,
 $$
 F_1^1\ge \frac{\b  |T-2t|}{t^2(T-t)^2}[c(T-2\d_0)-2R_1-\l^{-1})]= \frac{\b  |T-2t|}{t^2(T-t)^2}(cT-2R_1-2c\d_0-\l^{-1}),
  $$
which, via the second inequality in (\ref{2z}), is positive provided
that both $\d_0$ and $\l^{-1}$ are sufficiently small. Similarly,
for $t\in I_0$, $F_2^0$ is positive provided that $\d_0$ is
sufficiently small. Further, for $t\in I_0$,
 $$
  \ba{ll}
  G^1&\ds\ge \ds\frac{\b
 |T-2t|}{t^2(T-t)^2}\Big\{4c|T-2\d_0|\[c^2(\d_0-T/2)^2-R_1^2\]+O(\l^{-1})\Big\}\\
 \noalign{\ss}
 &\ge\ds \frac{\b
 |T-2t|}{t^2(T-t)^2}\Big\{4c|T-2\d_0|\[c^2T^2/4-R_1^2+c^2\d_0(\d_0-T)\]+O(\l^{-1})\Big\},
 \ea
 $$
which, via the second inequality in (\ref{2z}), is positive provided
that both $\d_0$ and $\l^{-1}$ are sufficiently small.

By (\ref{4.16})--(\ref{4.17}), we see that $F_1=F_1^0+F_1^1$,
$F_2=F_2^0+F_2^1$ and $G=G^0+G^1$. Noting the positivity of $F_1^0$,
$F_2^0$ and $G^0$, by the above argument, we see that $F_1$, $F_2$
and $G$ are positive for $t\in I_0$. For $t\in (0,T)\setminus I_0$,
noting again the positivity of $F_1^0$, $F_2^0$ and $G^0$, one can
choose $\b>0$ sufficiently small such that $F_1^1$, $F_2^1$ and
$G^1$ are very small so that $F_1$, $F_2$ and $G$ are positive.
Hence (\ref{4.14})--(\ref{4.15}) yield the desired (\ref{4.11}).
This completes the proof of Theorem \ref{4t1}.\endpf

\section{Proof of Theorems \ref{t1}-\ref{t2}}\label{s4}

We are now in a position to prove Theorems \ref{t1}-\ref{t2}.

\ms

{\it Proof of Theorem \ref{t2}}. The key idea is to apply Theorem
\ref{4t1}. Integrating both sides of (\ref{4.11}) (with $u$ replaced
by $y$, and $v=\th y$), using integration by parts, and recalling
that $\Th(t)$ decays exponentially to $0$ as $t\to0$ or $t\to T$,
noting that $v|_\Si=0$ (and hence $\n v=\frac{\pa v}{\pa\nu}\nu$ on
$\Si$), we arrive at
 \bel{4w11}
 \ba{ll}
 \ds\mathbb{E}\int_Q\left[c_0\l\Th \th^2(y_t^2+|\n y|^2+\l^2
y^2)+\Th\(-2\ell_tv_t+2\n\ell\cd\n v+\psi v\)^2\right]dxdt\\
 \noalign{\ss}
 \displaystyle
\le \mathbb{E}\int_Q\Th\th(-2\ell_tv_t+2\n\ell\cd\n v+\psi v)(dy_t-\D ydt)dx-\mathbb{E}\int_Q\Th \th^2\ell_t(dy_t)^2dx\\
\noalign{\ss}
\q\displaystyle+\mathbb{E}\int_\Si\Th\frac{\pa\ell}{\pa\nu}\left|\frac{\pa
v}{\pa\nu}\right|^2d\G dt.
 \ea\ee

By the first equation of system (\ref{hh6.1}), we get
 \bel{5w}
 \ba{ll}
 \ds \mathbb{E}\int_Q\Th\th(-2\ell_tv_t+2\n\ell\cd\n v+\psi v)(dy_t-\D ydt)dx-\mathbb{E}\int_Q\Th
 \th^2\ell_t(dy_t)^2dx\\
  \ns\ds
 =\mathbb{E}\int_Q\Th\th(-2\ell_tv_t+2\n\ell\cd\n v+\psi v)(a_1y_t+\lan a_2,\n y\ran+a_3y+f)dxdt\\
  \ns\ds\q-\mathbb{E}\int_Q\Th
 \th^2\ell_t(a_4y+g)^2dxdt\\
  \ns\ds
 \le \mathbb{E}\int_Q\Th\(-2\ell_tv_t+2\n\ell\cd\n v+\psi v\)^2dxdt\\
  \ns\ds\q+C\left\{\mathbb{E}\int_Q\Th\th^2\[a_1y_t+\lan a_2,\n y\ran+a_3y+f\]^2dxdt+\l\mathbb{E}\int_Q\Th
 \th^2(a_4y+g)^2dxdt\right\}\\
  \ns\ds
 \le \mathbb{E}\int_Q\Th\(-2\ell_tv_t+2\n\ell\cd\n v+\psi v\)^2dxdt\\
  \ns\ds\q
 +C\left\{\mathbb{E}\int_Q\Th\th^2(f^2+\l g^2)dxdt\right.+ |a_1|_{L^{\i}_{\cF}(0,T;(L^\i(G))}^2\mathbb{E}\int_Q\Th\th^2y_t^2dxdt\\
  \ns
  \ds\q+\l
  \[\l|a_3|_{L^{\i}_{\cF}(0,T;L^n(G))}^2+|a_4|_{L^{\i}_{\cF}(0,T;(L^\i(G))}^2\]
  \mathbb{E}\int_Q\Th
 \th^2y^2dxdt\\
  \ns
  \ds\q+\left.\[ |a_2|_{L^{\i}_{\cF}(0,T;L^{\i}(G;\mathbb{R}^n))}^2+|a_3|_{L^{\i}_{\cF}(0,T;L^n(G))}^2\]\mathbb{E}\int_Q\Th\th^2|\n
y|^2dxdt\right\}.
 \ea
 \ee

On the other hand, recalling (\ref{boundary}), we have
 \bel{215w}
 \ba{ll}\ds
 \mathbb{E}\int_\Si\Th\frac{\pa\ell}{\pa\nu}\left|\frac{\pa
v}{\pa\nu}\right|^2d\G
dt=2\l\mathbb{E}\int_\Si\Th\th^2(x-x_0)\cd\nu(x)\left|\frac{\pa
y}{\pa\nu}\right|^2d\G dt\\
 \ns
\ds\le
2\l\mathbb{E}\int_{\Si_0}\Th\th^2(x-x_0)\cd\nu(x)\left|\frac{\pa
y}{\pa\nu}\right|^2d\G_0 dt\le
C\l\mathbb{E}\int_{\Si_0}\Th\th^2\left|\frac{\pa
y}{\pa\nu}\right|^2d\G_0 dt.
 \ea
 \ee

Finally, combining (\ref{4w11}), (\ref{5w}) and (\ref{215w}), we
conclude the desired estimate (\ref{12e22}). This completes the
proof of Theorem \ref{t2}.\endpf

\ms

{\it Proof of Theorem \ref{t1}}. The proof follows easily from
Theorem \ref{t2} and the usual energy estimate. We omit the details.
\endpf

\end{document}